\documentclass[10pt,lettersize,journal]{IEEEtran}
\usepackage{amsmath,amsfonts}
\usepackage{algorithmic}
\usepackage{algorithm}
\usepackage{array}
\usepackage{textcomp}
\usepackage{stfloats}
\usepackage{url}
\usepackage{verbatim}
\usepackage{graphicx}
\usepackage{cite}
\usepackage{tikz}
\usepackage{mathrsfs}
\usepackage{makecell}
\usepackage{multirow}
\usepackage{threeparttable}
\usepackage{booktabs}

\usepackage{subfigure}

\newcommand{\A}{\mathcal{A}}
\newcommand{\B}{\mathcal{B}}
\newcommand{\C}{\mathcal{C}}

\newcommand{\E}{\mathcal{E}}
\newcommand{\F}{\mathcal{F}}
\newcommand{\G}{\mathcal{G}}
\newcommand{\N}{\mathcal{N}}
\newcommand{\V}{\mathcal{V}}
\newcommand{\W}{\mathcal{W}}

\newcommand{\mH}{\mathcal{H}}

\newtheorem{lemma}{Lemma}[section]
\newtheorem{property}{Property}[section]
\newtheorem{theorem}{Theorem}[section]
\newtheorem{problem}{Problem}[section]
\newtheorem{remark}{Remark}[section]

\newtheorem{definition}{Definition}[section]

\usepackage{color}

\newcommand{\black}{\textcolor{black}}

\newcommand{\block}{\mathrm{block}}
\newcommand{\diag}{\mathrm{diag}}

\begin{document}

\graphicspath{{figures/}}

\title{Event-Triggered Resilient Consensus of Networked Euler-Lagrange Systems Under Byzantine Attacks}

\author{Yuliang Fu, Guanghui Wen, Dan Zhao, Wei Xing Zheng,~\IEEEmembership{Fellow,~IEEE},~and~Xiaolei Li
%
\thanks{{This work was supported in part by the National Natural Science Foundation of China through Grant Nos. 62325304, U22B2046, U24A20279, 62306071, in part by the Jiangsu Provincial Scientific Research Center of Applied Mathematics under Grant No. BK20233002, and in part by the SEU Innovation Capability Enhancement Plan for Doctoral Students under Grant No. CXJH$\_$SEU25240. \textit{(Corresponding authors: Guanghui Wen; Wei Xing Zheng.)}}}
\thanks{Y. Fu is with the School of Cyber Science and Engineering, Southeast University, Nanjing 211189, China (e-mail: yuliang$\_$fu@163.com).}
\thanks{G. Wen is with the School of Automation, Southeast University, Nanjing 211189, China (e-mail: ghwen@seu.edu.cn).}
\thanks{D. Zhao is with the Department of Systems Science, Southeast University, Nanjing 211189, China (e-mail: danzhao@seu.edu.cn).}
\thanks{W. X. Zheng is with the School of Computer, Data and Mathematical Sciences, Western Sydney University, Sydney, NSW 2751, Australia (e-mail: w.zheng@westernsydney.edu.au).}
\thanks{X. Li is with the School of Electrical Engineering, Yanshan University, Qinhuangdao 066004, China (e-mail: xiaolei@ysu.edu.cn).}
}



\maketitle

\begin{abstract}
The resilient consensus problem is investigated in this paper for a class of networked Euler-Lagrange systems with event-triggered communication in the presence of Byzantine attacks. One challenge that we face in addressing the considered problem is the inapplicability of existing resilient decision algorithms designed for one-dimensional multi-agent systems. This is because the networked Euler-Lagrange systems fall into the category of multi-dimensional multi-agent systems with coupling among state vector components.
To address this problem, we propose a new resilient decision algorithm. This algorithm constructs auxiliary variables related to the coordinative objectives for each normal agent, and transforms the considered resilient consensus problem into the consensus problem of the designed auxiliary variables.
Furthermore, to relax the constraints imposed on Byzantine agent behavior patterns within continuous-time scenarios, the event-triggered communication scheme is adopted. Finally, the effectiveness of the proposed algorithm is demonstrated through case studies.
\end{abstract}


\section{Introduction}

In the past decade, the coordination of multi-agent systems (MASs) has attracted much attention from various scientific communities, owing primarily to its wide range of practical applications \cite{bk1}. In this context, agents exchange their local information through communication networks to accomplish desired collective tasks \cite{bk2}. As one of the most basic coordinative tasks, consensus has garnered increasing attention in recent years \cite{con1}.

It is noteworthy that the openness of communication networks exposes MASs to various cyber attacks. The presence of cyber attacks makes MASs difficult to achieve desired coordinative objectives. Recently, quite a few results on how to achieve consensus of MASs under several kinds of cyber attacks have been established \cite{att1}. Classical attacks include deception attacks, denial-of-service (DoS) attacks, false-data-injection attacks, replay attacks, and others. Yet Byzantine attacks considered in this paper represent a worse case due to their full capability to manipulate the behaviors of agents to launch more general and difficult class of attacks.

To achieve the consensus of MASs in the presence of Byzantine attacks, some resilient control schemes have been developed. Existing techniques can be generally divided into two categories: the Mean-Subsequence Reduced (MSR)-based approach \cite{vecbook} and the Resilient Vector Consensus (RVC)-based approach \cite{vce3}. In \cite{MSRE3}, the resilient cooperative output regulation problem under Byzantine attacks was addressed. The resilient consensus problem of switched MASs composed of continuous-time and discrete-time subsystems was solved in \cite{MSRE3n}. Furthermore, \cite{MSRE4} and \cite{MSRE4n} respectively considered the resilient consensus problem of MASs with state constraints and differential privacy requirements. The majority of the aforementioned MSR-based algorithms pertain to MASs with each agent being described by a one-dimensional dynamical system. Note that the MSR-based approaches become intricate when dealing with multi-dimensional MASs exhibiting coupling among state vector components.

For the RVC-based approach, some interesting results have been given in \cite{vce1,vce2,vce3}. Nevertheless, as noted therein, the graph robustness required by such an approach increases with the dimensions of the agents' dynamics. Particularly, as indicated in \cite{vecbook}, the achievable robustness is limited for any graph with a fixed number of agents. This means that such a method is incapable of dealing with the resilient consensus problem of MAS when the dimensions of each agents' dynamics exceed a specific value. The central point was used in \cite{vce4} to calculate the security point to improve the resilience of the RVC-based approach. However, the required computational resources increase significantly along with the increase of the system dimension. In light of these considerations, this paper focuses on exploring resilient consensus of networked Euler-Lagrange (EL) systems, driven by their widespread practical applications (such as rigid spacecraft, planar elbow manipulators, and marine vessels \cite{Spong2004book}). Despite their broad utility, the potential multi-dimensionality and intricate coupling among state vector components make it challenging to extend the aforementioned algorithms to effectively address the resilient consensus problem of networked EL systems.

Furthermore, an additional challenge in achieving consensus in MASs is that the aforementioned algorithms cause changes in their communication topologies. As mentioned in \cite{vecbook}, the existence of a minimum dwell time for these topologies cannot even be ensured in continuous-time systems. To address this, the study on continuous-time systems, including \cite{vecbook,MSRE3,MSRE3n,MSRE4}, constrains the information transmitted by Byzantine agents to normal agents to be continuous, and does not allow Byzantine agents to sever connections with their out-neighbors, which actually limits the capabilities of attackers. Noticeably, the event-triggered (ET) communication scheme provides a potential solution to this problem due to the existence of the minimum triggering interval and the open-loop estimation mechanism.

Recently, there have been some reported results on the consensus problem of MASs under ET communication schemes \cite{eti2,eti3,eti4,eti5}.
For example, the work in \cite{eti3} proposed a velocity-free ET scheme to deal with the consensus problem of networked EL systems with communication delays. In \cite{eti2}, a dynamic ET scheme was designed to achieve consensus in networked EL systems with uncertain parameters. The fault-tolerant consensus problem of nonlinear MASs with ET communication was solved in \cite{eti4} by using the dynamic ET scheme. Additionally, an ET scheme was devised in \cite{eti5} based on the axis-angle attitude representation to solve the finite-time attitude consensus problem of MASs under switching topologies. Nevertheless, these results are generally inapplicable to the resilient consensus problem, as the presence of Byzantine attacks inhibits normal agents from accessing reliable global information and may even lead to different global information in different dimensions.

Motivated by the above discussions, this paper focuses on solving the resilient consensus algorithm for networked EL systems with state coupling under Byzantine attacks in a continuous-time scenario. The inherent state coupling in these systems, along with the potential for Byzantine agents to withhold transmissions or transmit discontinuous information to their out-neighbors, makes existing algorithms unsuitable.

The main contributions of this paper are threefold. First, a new resilient decision algorithm is proposed for multi-dimensional MASs with state coupling. This algorithm transforms the considered problem into the consensus problem of the designed auxiliary variables, whose dimensions are independent of one another, thereby simplifying the analysis in scenarios with state coupling and reducing the computational resources required. The reduced execution frequency of the proposed algorithm resulting from the ET communication scheme further conserves the computational resources. These technical innovations make the proposed algorithm more feasible for multi-dimensional MASs with state coupling than those in \cite{MSRE3,MSRE3n,MSRE4,MSRE4n,vce1,vce2,vce3}, which either fail to handle state coupling, impose strict constraints on attacks, or exhibit high demands on graph robustness and computational resources.
Second, by utilizing the ET communication scheme, the communication topologies used in the consensus problem of auxiliary variables are ensured to be connected, and their Laplacian matrices are guaranteed to be piecewise continuous, regardless of the behavior of Byzantine agents. Compared to \cite{vecbook,MSRE3,MSRE3n,MSRE4,MSRE4n,vce1,vce2,vce3}, the proposed algorithm neither requires the information transmitted by Byzantine agents to lie on a continuous trajectory, nor requires the agents to maintain continuous connections with their out-neighbors, thus relaxing the constraints on Byzantine attacks.
Finally, a fully distributed resilient consensus algorithm with ET communication under a digraph is devised, ensuring the system security even in scenarios of heterogeneous topologies with different dimensions, further enhancing the technical feasibility of the proposed algorithm.

\emph{Notations}: Let $\mathbb{R}^n$ be the set of real column vectors with dimension $n$ and $\mathbb{R}^{n \times m}$ be the set of $n\times m$ real matrices. $\mathbb{R}$ and $\mathbb{N}$ represent set of real numbers and non-negative integers, respectively. $\mathbf{0}_n\in \mathbb{R}^n$ and $\mathbf{1}_n\in \mathbb{R}^n$ are $n$-dimensional column vectors with all elements being 0 and 1, respectively. And $I_n \in \mathbb{R}^{n \times n}$ is the identity matrix with dimension $n$. Let $\otimes$ be the Kronecker product, $\circ$ be the Hadamard product and $\mathcal{O}$ be the Bachmann-Landau notation. $|| \cdot ||$ is the Euclidean norm of a vector or the Frobenius norm of a matrix. $X^T$ is the transpose of $X$, where $X$ can be a vector or a matrix. Denote by $\block \, \diag(X_1,X_2,\cdots,X_N)$ a block diagonal matrix with the entries $X_i$, where $i=1,\cdots ,N$, and $X_i$ can be a scalar or a matrix. When $X$ is a scalar, define $\lceil X \rceil = \min\{m \in \mathbb{N^+}\mid m \ge X\}$, where $\mathbb{N^+}$ is the set of positive integers. For a function $g(t):[0,+\infty) \to \mathbb{R}$, define its upper Dini derivative as $D^+ g(t)=\lim \sup\nolimits_{\delta \to 0^+} \frac{g(t+\delta)-g(t)}{\delta}$. For any two sets $\A$ and $\B$, $\A \bigcap \B$, $\A \bigcup \B$, and $\A \backslash \B$ represent their intersection set, union set, and difference set, respectively. Additionally, $|\A|$ is the number of elements in a set $\A$.

\section{Preliminaries and Problem Statement}\label{sec_pf}

\subsection{Graph Theory}

This paper considers a networked EL system comprising $N$ agents. In such a system, the information exchange among agents is assumed to be represented by a digraph, denoted as $\G = \{\V,\E\}$, where $\V = \{v_i \mid i=1,\cdots, N\}$ is the vertex set and $\E \subset \V \times \V$ is the edge set. The edge $(v_i,v_j) \in \E$ signifies that agent $i$ can receive the information transmitted by agent $j$, and the in-neighbor set of agent $i$ is defined as $\N_i=\{j \mid (v_i,v_j) \in \E\}$. Define $A=[a_{ij}] \in \mathbb{R}^{N \times N}$ as the adjacency matrix, where $a_{ij}=1$ if $j \in \N_i$ and $a_{ij}=0$ otherwise. Similarly, $D=[d_{ij}] \in \mathbb{R}^{N \times N}$ is defined as the in-degree matrix, where $d_{ij}=\sum\nolimits_{j \in \N_i}{a_{ij}}$ if $i=j$ and $d_{ij}= 0$ otherwise. Accordingly, the Laplacian matrix $L$ of digraph $\G$ can be defined as $L=D-A$.

In general, the robustness of a communication topology is crucial to guarantee the achievement of consensus in MASs under Byzantine attacks. Accordingly, the following definitions are introduced to characterize the communication topology robustness \cite{MSRE3}.

\begin{definition}[\textbf{$r$-reachable}] \label{dfreachable}
Consider a nonempty subset $S_0$ of the vertex set $\V$ in a digraph $\G=\{\E,\V\}$. For an $r \in \mathbb{N}$, if there exists an agent $i\in S_0$ such that $|\N_i\backslash S_0| \ge r$, then $S_0$ is referred to as an $r$-reachable set.
\end{definition}

\begin{definition}[\textbf{$r$-robust}] \label{dfrobust}
Suppose that $S_1$ and $S_2$ are any pair of nonempty disjoint subsets of $\V$. The digraph $\G$ is said to be $r$-robust if at least one of $S_1$ and $S_2$ is $r$-reachable.
\end{definition}


\begin{figure}[!t]
\centering
\subfigure[A 2-robust digraph]{
    \includegraphics[width=0.30\linewidth]{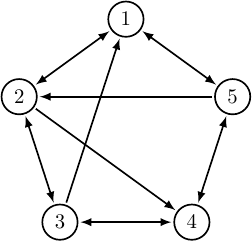}}~
\subfigure[Under 1-total Byzantine attack]{
    \includegraphics[width=0.30\linewidth]{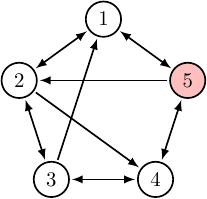}}~
\subfigure[Under 1-local Byzantine attack]{
    \includegraphics[width=0.30\linewidth]{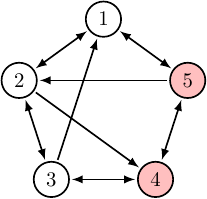}}
\caption{An example of 2-robust digraph with 5 agents, 
where the red agents are Byzantine agents.}
\label{2robusttgraph}
\end{figure}

To help understand these two definitions, Fig.~\ref{2robusttgraph}(a) presents a 2-robust digraph consisting of 5 agents. In this figure, let us choose $S_0{=}\{1,5\}$. Since $\N_1{=}\{2,3,5\}$ and $\N_5{=}\{1,4\}$, one can get $|\N_1\backslash S_0|{=}|\{2,3\}|{=}2$ and $|\N_5\backslash S_0|{=}1$. Therefore, according to Definition~\ref{dfreachable}, we can say that $S_0$ is a 2-reachable set because $|\N_1\backslash S_0|{=}2$. Additionally, since $|\N_1\backslash S_0|{>}1$, $S_0$ is also 1-reachable.

Moreover, based on the rule in Definition~\ref{dfrobust} that the sets $S_1$ and $S_2$ are nonempty and disjoint, we observe that when one of these sets contains 3 or 4 agents, the other must contain 1 or 2 agents. Therefore, we only need to analyze the reachability of all possible subsets of $\V$ that contain 1 or 2 agents. From Fig.~\ref{2robusttgraph}(a), it is evident that when either $S_1$ or $S_2$ contains only one agent, it is always 2-reachable, as each agent has no fewer than two in-neighbors. When $S_1$ or $S_2$ contains two agents, if one of the agents is not an in-neighbor of the other, then this conclusion is also straightforward. For the remaining cases, where the set is one of $\{1,2\}$, $\{2,3\}$, $\{3,4\}$, $\{4,5\}$, or $\{5,1\}$, it is also 2-reachable according to Definition~\ref{dfreachable}. Therefore, the topology in Fig.~\ref{2robusttgraph}(a) is 2-robust.

\subsection{Euler-Lagrange Systems}

The dynamics of each agent $i \in \V$ can be described by the following EL equation:
\begin{equation}\label{eq_ELequation}
{M_i}({q_i}){\ddot q_i} + {C_i}({q_i},{\dot q_i}){\dot q_i} + {g_i}({q_i}) = {\tau_i},
\vspace{-0.5ex}
\end{equation}
where $q_i \in \mathbb{R}^n$ and $\tau_i \in \mathbb{R}^n$ are the generalized position and the control input of agent $i$, respectively. ${M_i}({q_i}) \in \mathbb{R}^{n \times n}$ is the inertia matrix, ${C_i}({q_i},{\dot q_i})\in \mathbb{R}^{n \times n}$ is the Coriolis and centripetal matrix, and $g_i \in \mathbb{R}^n$ is the gravity vector of agent $i$. Some useful properties of EL dynamics are given below.

\begin{property}\label{pM}
The inertia matrix ${M_i}({q_i})$ is a bounded, symmetric positive-definite matrix.
\end{property}

\begin{property}\label{p1}
For any vector $p \in \mathbb{R}^n$, there holds ${p^T} ({\dot M_i}({q_i}) - 2{C_i}({q_i},{\dot q_i})) p =0$.
\end{property}

\begin{property}\label{p2}
Let $\varphi_i \in \mathbb{R}^l$ be a vector composed of system parameters, with the elements of $\varphi_i$ being all unknown. Let $\Omega_i( q_i, \dot q_i, x_i, y_i) \in \mathbb{R}^{n \times l}$ be a regression matrix, which can be obtained directly. Then, for any measurable $x_i \in \mathbb{R}^n$ and $y_i \in \mathbb{R}^n$, \eqref{eq_ELequation} can be rewritten as
\begin{equation}\label{eq_ELlinear}
{M_i}({q_i}){x_i} + {C_i}({q_i},{\dot q_i})y_i + {g_i}({q_i}) = \Omega_i \varphi_i.
\end{equation}
\end{property}

\begin{remark}
EL dynamics are widely used to model practical systems, including rigid spacecrafts, planar elbow manipulators, and marine vessels \cite{Spong2004book}. However, the multi-dimensionality and complex nonlinearity of EL dynamics pose challenges for extending the existing resilient consensus algorithms to such systems. As a consequence, addressing the resilient consensus problem in networked EL systems becomes meaningful.
\end{remark}

\subsection{Attack Model}

Byzantine attacks aim at degrading the algorithm's performance by preventing the achievement of global objectives. In general, to guarantee the achievement of resilient consensus of normal agents in a worst-case network environment, attackers are assumed to be omniscient and can have full capability to manipulate specific agents, referred to as Byzantine agents \cite{arxiv1}. The possible behavior pattern of Byzantine agents can be summarized as follows.
\begin{itemize}
  \item Byzantine agents do not adhere to the algorithms designed for normal agents.
  \item Byzantine agents can arbitrarily update their states and information transmitted to each of their out-neighbors according to the schemes planned by the attackers to disrupt the consensus of normal agents, perhaps by colluding with each other to do so.
\end{itemize}

The aforementioned model finds widespread adoption in modeling Byzantine agents, as observed in \cite{MSRE4n,MSRE3,MSRE3n,MSRE4}. However, launching such attacks requires significant resource consumption, constraining the maximum number of Byzantine agents that attackers with limited resources can manipulate. Without this constraint, this problem could become overly pessimistic and inconsequential \cite{MSRE3}. To model the number and distribution of Byzantine agents, we define $\F$ as the set of Byzantine agents and $\mH = \V \backslash \F$ as the set of normal agents. Then we introduce the following two models.

\begin{definition}[\textbf{$f$-total and $f$-local attack} \cite{MSRE3}]\label{dffattack}
For an $f \in \mathbb{N}$, an MAS is said to be under an $f$-total attack if there exist at most $f$ Byzantine agents in it, i.e., $|\F| \le f$. Further, if $|\N_i \bigcap \F| \le f$ holds for any normal agent $i$, then the system is said to be under an $f$-local attack.
\end{definition}

The above two models were initially proposed in \cite{fmodel} and have since been widely adopted in \cite{MSRE4n,MSRE3,MSRE3n,MSRE4}. In addition, in this paper, we have relaxed the assumption in \cite{vecbook,MSRE3,MSRE3n,MSRE4}, allowing that the information transmitted by Byzantine agents to their out-neighbors in continuous-time systems can be discontinuous and Byzantine agents can arbitrarily sever connections with their out-neighbors. This enables the proposed algorithm to withstand more aggressive Byzantine attacks.

\begin{remark}
It is noteworthy that attacks satisfying the $f$-total attack model always satisfy the $f$-local attack model, but not vice versa. This is because, even though $|\N_i \bigcap \F| \le f$ holds for any normal agent $i$, when $|\N_i \bigcap \N_j|\le f$, there may be cases where \black{$|(\N_i \bigcap \F) \bigcup (\N_j \bigcap \F)| \ge f$}. Therefore, there may exist more than $f$ Byzantine agents within the entire network. In other words, the $f$-total attack model is a specific case of the $f$-local attack model. Based on the aforementioned discussion, this paper concentrates on the resilient consensus problem of networked EL systems in the presence of $f$-local Byzantine attacks. Additionally, to make this conclusion evident, Fig.~\ref{2robusttgraph}(b)--(c) present an example of a 2-robust digraph subject to 1-total Byzantine attack and 1-local Byzantine attack.
\end{remark}


\subsection{Problem Formulation}

This paper is devoted to developing a fully distributed event-triggered resilient consensus algorithm for networked EL systems under Byzantine attacks. Unlike the traditional consensus problem in \cite{eti2,eti3,eti4,eti5}, the resilient consensus problem requires consensus among normal agents rather than all agents, as described in \cite{vecbook,vce3,MSRE3,MSRE3n,MSRE4,MSRE4n,vce1,vce2,vce4}, since Byzantine agents do not follow the preset algorithm. The control objectives can be summarized as follows.

\begin{problem}\label{df_obj1}
To achieve the resilient consensus of networked EL systems, how to design the resilient consensus algorithm for each normal agent so as to guarantee that
\vspace{-0.5ex}
\begin{equation}\label{eq_objcon}
\mathop {\lim }\limits_{t \to \infty } (q_j(t) - q_i(t)) = \mathbf{0}_n, \;
\mathop {\lim }\limits_{t \to \infty } (\dot q_j(t) - \dot q_i(t)) = \mathbf{0}_n,
\vspace{-0.5ex}
\end{equation}
as well as to ensure that $q_i(t)$ and $\dot q_i(t)$ remain bounded for any $t$, where $i,j \in \mH$. Additionally, owing to the utilization of the ET communication scheme, it is imperative to guarantee elimination of the Zeno behavior.
\end{problem}

\section{Main Results}

In this section, a new resilient decision algorithm is first proposed for networked EL systems, addressing their resilient consensus problem in the presence of Byzantine attacks. Meanwhile, a fully distributed consensus algorithm is designed for each agent to complement the proposed resilient consensus algorithm. Additionally, the ET communication scheme is also adopted to reduce the communication resource consumption.

\subsection{Design of Resilient Consensus Algorithm}

In this subsection, a resilient consensus algorithm for networked EL systems under ET communication is proposed. This algorithm includes a new resilient decision algorithm, named as the auxiliary-variable-based resilient decision (AVBRD) algorithm, and a fully distributed consensus algorithm. We begin by detailing the proposed AVBRD algorithm.

The AVBRD algorithm guarantees the achievement of consensus in networked EL systems under Byzantine attacks by introducing an auxiliary variable $\hat {\W}_i \in \mathbb{R}^{n}$ related to the preset control objectives, transforming the resilient consensus problem of networked EL systems into the consensus problem of $\hat {\W}_i$. Given that Byzantine agents aim to prevent the achievement of preset control objectives in networked EL systems, the auxiliary variable should be designed to gradually achieve consensus as the global objectives tend towards being achieved. For example, for the consensus algorithm with absolute velocity damping of network EL systems, $\W_i$ can be constructed based on $y_i$ or $x_i$, which are defined in \cite{MJnew}; for the static scalar consensus problem, $\hat {\W}_i$ can be constructed based on $k_iq_i$ with $k_{i}$ being the scaled scalar.
After constructing the auxiliary variable, agent $i$ calculates the mean of the ($f$+1)-th and the ($|\N_i|{-}f$)-th smallest or largest values of $\hat {\W}^l_j(t)$ in each dimension $l$ for all $j \in \N_i$, where $\hat {\W}^l_j(t)$ is the $l$-th dimension of $\hat {\W}_j(t)$. Then, it uses the difference between the obtained value and $\hat {\W}_i$ to drive the control input. The global objectives are thus transformed into the consensus problem of the auxiliary variables. The algorithmic details are presented in Algorithm~\ref{alg-1}.

\begin{algorithm}[!t]
\caption{The AVBRD Algorithm {from the Viewpoint of Normal Agent $i$}}\label{alg-1}
\begin{algorithmic}[1]
\STATE 
Agent $i$ calculates $ \hat {\W}_j(t)$ for each of its in-neighbors, based on the information it stored.
\STATE 
For $j_k \in \N_i$ with $k=\{1,2,\cdots,|\N_i|\}$, agent $i$ constructs matrix $[ \hat {\W}_{j_1}(t), \hat {\W}_{j_2}(t), \cdots, \hat {\W}_{j_|\N_i|}(t)] \in \mathbb{R}^{n\times |\N_i| }$, and sorts each dimension on the same rule to obtain $ \hat {\boldsymbol{{\W}}}_i (t)$.
\STATE 
Agent $i$ calculates the mean of $ \hat {\W}(t)$ of the ($f$+1)-th and the ($|\N_i|{-}f$)-th columns on each of its corresponding row independently, and the result is given as $\bar{\W}_i(t)$.
\end{algorithmic}
\end{algorithm}

Obviously, the introduction of the resilient decision algorithm may lead to changes in the global information. Therefore, it is necessary to design a fully distributed consensus algorithm without using any global information for each agent. However, the complex nonlinearities present in networked EL systems make the development of such an algorithm challenging with ET communication under the digraph. To combat this challenge, we adopt the output regulation technique. The designed algorithm comprises a fully distributed ET observer and a distributed dynamic control law. Next, we introduce the following observer to generate a common reference trajectory for normal agents:
\begin{equation}\label{eq_obs}
\dot \eta_i(t)= S \eta_i(t) - \mu_1 (\hat {\eta}_i(t) - \bar{\eta}_i(t)),
\end{equation}
where $\mu_1 >0$ is a constant gain, and $S \in \mathbb{R}^{n\times n}$ is a semi-simple system matrix whose eigenvalues all have zero real parts, implying the state coupling in the observer. $\bar{\eta}_i(t)$ is derived based on the value of $\bar {\W}_i(t)$, which has been defined in Algorithm~\ref{alg-1}. Their specific mapping depends on that between $\W_i$ and $\eta_i(t)$. $\hat {\W}_i(t)$ and $\hat {\W}_j(t)$ are the auxiliary variable deigned based on $\hat \eta_i(t)$ and $\hat \eta_j(t)$, where $\hat \eta_i(t)$ and $\hat \eta_j(t)$ are the open-loop estimates of $ \eta_i(t)$ and $ \eta_j(t)$, respectively, and they are defined as
\begin{equation}\label{eq_hateta}
\hat \eta_i(t) = e^{S(t-t^i_k)} \eta_i(t^i_k), \ \ \hat \eta_j(t) = e^{S(t-t^j_m)}\eta_j(t^j_m),
\end{equation}
in which, $t^i_k$ is the $k$-th triggering instant of agent $i$, and $t^j_m$ is the latest triggering instant of agent $j$ before $t$.

It is worth mentioning that the above closed-loop estimation can be completed locally by agent $i$, as agent $j$ only needs to transmit information to agent $i$ at its triggering instants. In this paper, obtaining the state $q_j(t)$ of its in-neighbor $j$ in real-time is challenging for normal agents. Fortunately, by applying \eqref{eq_event} and \eqref{eq_controller} to be given, achieving the control objective \eqref{eq_objcon} will also ensure that $\mathop {\lim }\nolimits_{t \to \infty } (\eta_j(t) - \eta_i(t)) = \mathbf{0}_n, \forall i,j \in \mH$. Therefore, we construct $\W_i(t)= e^{-St} \eta_i(t)$ and $\hat {\W}_i(t)= {\W}_i(t_k^i)$, which can be re-expressed as $\hat {\W}_i(t)= e^{-St}{\hat \eta}_i(t)$. Accordingly, we devise $\bar \eta_i(t)=e^{St} \bar {\W}_i(t)$. Inspired by \cite{Liutac}, the following fully distributed ET scheme is designed:
\begin{equation}\label{eq_event}
t_{k + 1}^{i} = \inf \Big\{ {t > t_k^{i}\Big| \|e_{\eta i}(t)\| - \frac{\alpha_{1}}{(t-t_0+\alpha_2)^{\alpha_3}} \ge 0} \Big\},
\end{equation}
where $\alpha_{1}>0$, $\alpha_{2}>1$ and $\alpha_{3}>1$ are all constant gains, $t_0$ is the initial instant, and $e_{\eta i}(t)= \hat \eta_i(t)-\eta_i(t)$ is the error variable. The distributed dynamic control law is devised as
\begin{equation}\label{eq_controller}
\left\{\begin{array}{l}
\tau_i = -k_i s_i + \Omega_i( q_i, \dot q_i, v_i, \dot v_i) \hat \varphi_i,\\
 v_i = S \eta_i -\mu_{2} (q_i -\eta_i) ,\\
\dot {\hat \varphi}_i = -F_i \Omega_i^T s_i,
\end{array} \right.
\end{equation}
where $v_i \in \mathbb{R}^n$ is defined as a reference velocity, $s_i = \dot q_i - v_i$ is defined as an auxiliary variable, $F_i>0$ and $\mu_2>0$ are both constant control gains, $\hat \varphi_i$ is the estimate of $\varphi_i$, and $\Omega_i$ and $\varphi_i$ have been given in Property \ref{p2}. The proposed resilient consensus algorithm is summarized as Algorithm~\ref{alg-2}.

\begin{algorithm}[!t]
\caption{The Proposed Resilient Consensus Algorithm from the Viewpoint of Normal Agent $i$}\label{alg-2}
\begin{algorithmic}[1]
\STATE 
Agent $i$ checks whether any of its in-neighbors have transmitted their state information to it. If so, update the storage for those in-neighbors that have a time interval since their last transmit.
\STATE 
Agent $i$ calculates $\hat \eta_i(t)$ and $\hat \W_i(t)$ for itself. Then, if its storage is updated, agent $i$ executes the AVBRD algorithm.
\STATE 
Agent $i$ measures its observer state $\eta_i(t)$, as well as its position state $q_i(t)$ and velocity state $\dot q_i(t)$, and updates $\hat \varphi_i$. Based on \eqref{eq_controller}, Property~\ref{p2} and the above measurements, agent $i$ calculates the new values of $v_i$, $s_i$, $\Omega_i$, and $\dot {\hat \varphi}_i$.
\STATE 
Agent $i$ calculates the values of $\bar \eta_i(t)$, $\dot \eta_i(t)$, and $\tau_i(t)$ according to the above measurements and calculation results, as well as \eqref{eq_obs}, \eqref{eq_controller} and the information stored by agent $i$.
\STATE 
Agent $i$ checks whether $t$ is a triggering instant according to \eqref{eq_event}. If so, agent $i$ transmits its observer state $\eta_i(t)$ and the triggering instant $t$ to all its out-neighbors.
\end{algorithmic}
\end{algorithm}

\begin{remark}\label{rmk31}
It is noteworthy that, for multi-dimensional MASs, when both $n$ and $|\N_i|$ are positive integers not less than 2, both Algorithm~\ref{alg-1} and the MSR-based algorithm have the time complexity of $\mathcal{O}(n|\N_i|\log_2|\N_i|)$. However, due to the introduction of the ET communication scheme and the design of $\hat {\W}_i(t)= {\W}_i(t_k^i)$, Algorithm~\ref{alg-1} does not need to be executed at all instants, which significantly reduces the computational resources it requires. Moreover, although the auxiliary variables are required to be computed, its time complexity, as well as that of the subtraction operations in the MSR-based algorithm, is all $\mathcal{O}(m)$, where $m$ is the number of elements to be computed. Additionally, the auxiliary variables corresponding to the new information are computed only at the instant when they are received by agent $i$, which further reduces the resource consumption of the proposed algorithm.
\end{remark}

\begin{remark}\label{rmk32}
Note that the RVC-based algorithm can be applied only when $|\N_i|\ge (n+1)f+1$. In this case, according to \cite{vce3}, the time complexity of this algorithm is $\mathcal{O}(n|\N_i|\log_2|\N_i|)+\mathcal{O}((r_1r_2)^3)$, where $r_1=nf+1$ and $r_2=\binom{(n+1)f+1}{f}$. The term $\mathcal{O}(n|\N_i|\log_2|\N_i|)$ arises from the sorting process of the points, while the other term attributes to the computation of the convex hull. This formulation is due to the unclear specific values of $|\N_i|$, $f$ and $n$. Clearly, the proposed algorithm requires less computational resources.
\end{remark}

\begin{remark}
Note that, in the existing MSR-based algorithms for continuous-time systems \cite{vecbook,MSRE3,MSRE3n,MSRE4}, Byzantine attacks are imposed with strong continuity assumption, i.e., the information transmitted by a Byzantine agent to its any out-neighbor must lie on a continuous trajectory, and Byzantine agents are prohibited from severing their connections with their out-neighbors. The proposed algorithm relaxes these constraints. Additionally, the RVC-based algorithms like \cite{vce1, vce2, vce3} require the communication topology to be (($n$+1)$f$+1)-robust. Following \cite[Property~5.19]{vecbook}, no graph with $N$ agents can exhibit ($\lceil \frac{N}{2} \rceil$+1)-robustness. This implies that the RVC-based algorithms become ineffective when $n>\lceil \frac{N}{2f} \rceil -1$.
In summary, the proposed algorithm exhibits broader applicability than the existing algorithms.
\end{remark}

\begin{remark}
Note that although \cite{Wangtac} proposed a similar observer under continuous communication, global information plays a crucial role in the convergence of $\eta_i(t)$ therein. In addition, though the ET scheme in \eqref{eq_event} adopts a similar form as that in \cite{Liutac}, the potential heterogeneity of topology in different dimensions in \eqref{eq_lmobseq} caused by Algorithm~\ref{alg-1} introduces additional difficulties for the theoretical analysis presented in this paper. Therefore, the aforementioned results in \cite{Wangtac,Liutac} cannot be easily extended to the scenarios considered in this paper.
\end{remark}

To help develop the analysis process, we introduce the following lemmas.

\begin{lemma}[\!\!\!\cite{lmfrom1}]\label{lm_seclm1}
If the communication graph $\G$ with $N$ agents contains a spanning tree, then the spectrum of its Laplacian matrix $L$ consists of a simple zero eigenvalue and $N-1$ eigenvalues with positive real parts.
\end{lemma}

\begin{lemma}[\!\!\!\cite{Liutac}]\label{lm_seclm2}
Consider that $g(t):[t_0,+\infty) \to \mathbb{R}$ is a continuous, nonnegative and bounded function, with $\lim \nolimits_{t \to \infty} g(t)=0$. Then, for any $\omega_1>0$, there exists a positive constant $\omega_2$ such that $\int \nolimits_{\tau}^{t} g(s) ds \le \omega_1(t-\tau) +\omega_2$, where $t\ge\tau\ge t_0$.
\end{lemma}

\begin{lemma}\label{lm_neighbours}
If the communication graph $\G$ is an $r$-robust graph with $r>1$, then $|\N_i| \ge r$ for any $i \in \V$.
\end{lemma}

\begin{IEEEproof}[\hspace{-1.0em}Proof]
{See the Appendix.}
\end{IEEEproof}

\begin{lemma}\label{lm_obsever2lm}
In a (2$f$+1)-robust graph $\G$ under $f$-local Byzantine attacks, for any of its normal agent $i$ with dimension $l$, there holds that $\bar{\W}^l_i(t)$ lies within the convex hull formed by the auxiliary variables $\hat{\W}^l_j(t)$ from at least $|\N_i|-2f$ normal agents $j$ among its in-neighbors.
\end{lemma}

\begin{IEEEproof}[\hspace{-1.0em}Proof]
{See the Appendix.}
\end{IEEEproof}

\begin{lemma}\label{lm_obsever}
Assume that there exist Byzantine agents. Suppose that for any agent $i$, the observer \eqref{eq_obs} is written in the form as
\begin{equation}\label{eq_lmobseq}
\dot{\W}_i(t)= - \mu_1 \sum \nolimits_{j=1}^{N} \tilde a_{ij}(t)\circ(\hat {\W}_i(t) - \hat{\W}_j(t)),
\end{equation}
where $\tilde a_{ij}(t)=[\tilde a_{ij}^1(t),\tilde a_{ij}^2(t),\cdots,\tilde a_{ij}^n(t)]^T$, and $\tilde a_{ij}^l(t)$ is the element of the adjacency matrix of the graph $\tilde \G^l(t)$, which contains at least one directed spanning tree and whose adjacency matrix is piecewise continuous. Moreover, the constant gains mentioned in \eqref{eq_obs} and \eqref{eq_event} are chosen as $\mu_1>0$, $\alpha_1>0$, $\alpha_2>1$\black{,} and $\alpha_3>1$. Then, the observer \eqref{eq_obs} and the ET scheme \eqref{eq_event} guarantee that the observer states of each agent asymptotically converge to the same value. Additionally, the Zeno behavior is eradicated.
\end{lemma}

\begin{IEEEproof}[\hspace{-1.0em}Proof]
	See the Appendix.
\end{IEEEproof}

\subsection{Analysis of Resilient Consensus}

The main result of this work is presented as follows.

\begin{theorem}\label{them1}
Consider a networked EL system consisting of $N$ agents, where each agent is modeled by the EL equation \eqref{eq_ELequation}. If each agent obtains the initial state of all its in-neighbors and the communication topology is (2$f$+1)-robust, then the event-triggered resilient consensus problem under $f$-local Byzantine attacks is solved by the proposed resilient consensus algorithm (see Algorithm~\ref{alg-2}), which consists of the observer \eqref{eq_obs}, the ET scheme \eqref{eq_event}, the dynamic control law \eqref{eq_controller} and the AVBRD algorithm. Moreover, the choice of the constant gains mentioned in \eqref{eq_obs} and \eqref{eq_event} is the same as in Lemma~\ref{lm_obsever}, and the parameters mentioned in \eqref{eq_controller} and Algorithm~\ref{alg-1} are chosen as $\W_i=e^{-St}\eta_i$, $\mu_{2} > 0$, $k_i> 0$, $F_i> 0$, and $k_i \mu_{2} > \frac{1}{2}$.
\end{theorem}

\begin{IEEEproof}[\hspace{-1.0em}Proof]
By combining with \eqref{eq_obs}, the derivative of $\W_i$ can be obtained as
\begin{equation}\label{eq_dW}
\dot \W_i(t) = -Se^{-St}\eta_i + e^{-St} \dot \eta_i = - \mu_1 (\hat {\W}_i(t) - \bar{\W}_i(t)).
\end{equation}

To apply the conclusion from Lemma~\ref{lm_obsever}, we first attempt to rewrite \eqref{eq_dW}. From this, we have $\dot \W_i^l(t) = - \mu_1 (\hat {\W}^l_i(t) - \bar{\W}^l_i(t))$, where $\bar {\W}^l_i(t)$ is the $l$-th dimension of $\bar {\W}_i(t)$, and $\W^l_i(t)$ and $\hat {\W}^l_i(t)$ have been defined in the proof of Lemma~\ref{lm_obsever}. Recalling Lemma~\ref{lm_obsever2lm}, it follows from Algorithm~\ref{alg-1} that there holds $\bar{\W}^l_i(t) = \sum_{j\in \mH} \tilde a^l_{ij}(t) \hat{\W}^l_{j}(t)$ at any dimension $l$ and each instant $t$, where $ \sum_{j\in \mH} \tilde a^l_{ij}(t)=1 $ and $ \tilde a^l_{ij}(t) \in [0,1]$. Here, we emphasize that, even if Byzantine agents disconnect from their out-neighbors during specific periods, this result still holds. This is because, for each normal agent $i$, Algorithm~\ref{alg-2} ensures that $\hat {\boldsymbol{{\W}}}_i (t)$ maintains the same number of columns after the initial instant. Thus, \eqref{eq_dW} can be expressed in the form of \eqref{eq_lmobseq}, and the in-neighbors of a normal agent in any $\tilde \G^l(t)$ exclude Byzantine agents.

Next, let us define $\tilde \G^l_{\mH}(t)$ as the subgraph of $\tilde \G^l(t)$ that only contains all the normal agents in $\tilde \G^l(t)$ and the edges among them. We will now provide a proof that the there always exists a directed spanning tree among normal agents in $\tilde \G^l_{\mH}(t)$. Suppose that $\G$ is a (2$f$+1)-robust graph. For any pair of nonempty disjoint vertex sets $S_1$ and $S_2$ of the digraph $\G$, there exists at least one agent $i \in S_p$, $p=\{1,2\}$ such that $|\N_i \bigcap (\V \backslash S_p) | \,{\ge}\, 2f+1$. Based on the aforementioned discussion, after using the AVBRD algorithm, an agent $i$ in $\tilde \G^l(t)$ retains its at least $|\N_i|-2f$ in-neighbors for any dimension $l$ at each instant. Therefore, for the mentioned agent $i$, there holds $|\tilde \N^l_i(t) \bigcap (\V \backslash S_p) | \,{\ge}\, {1}$ regardless of whether $f\,{=}\,0$, where $\tilde \N^l_i(t)$ is the in-neighbor set of agent $i$ in $\tilde \G^l(t)$. That is, at least one of {$S_1$ and $S_2$} is ${1}$-reachable, so $\tilde \G^l(t)$ is a $1$-robust graph, indicating that it contains a directed spanning tree \cite[Property~5.14]{vecbook}. Furthermore, note that in any $\tilde \G^l(t)$, Byzantine agents are not in-neighbors of normal agents, meaning that after removing all Byzantine agents, the directed spanning tree still exists. Moreover, owing to the open-loop estimation mechanism \eqref{eq_hateta}, such conclusions will always hold as long as each agent can obtain the initial states of all its in-neighbors.

Then, \eqref{eq_dW} can be written into
\begin{equation}\label{eq_dW2}
\dot \W_i(t) = - \mu_1 \sum \nolimits_{j\in \mH} \tilde a_{ij}(t) \circ (\hat {\W}_i(t) - \hat{\W}_j(t)),
\end{equation}
which is equal to \eqref{eq_lmobseq} with $\tilde a_{ij}(t)=\mathbf{0}_n, \forall j \notin \N_i\bigcap \mH$.
Following Lemma~\ref{lm_obsever}, when $\W_i(t_0)$ is bounded and Zeno behavior does not occur at $t_0$, normal agents can hold a minimum triggering interval. And, the first step of Algorithm~\ref{alg-2} ensures that even if Byzantine agents continuously update information to their out-neighbors, normal agents will not update their storage. Combined with the existence of a minimum triggering interval, this mechanism does not affect the information from normal agents. Altogether, these ensure that, when $\W_i(t_0)$ is bounded and Zeno behavior does not occur at $t_0$, the adjacency matrix of $\tilde \G^l_{\mH}(t)$ is piecewise continuous and row-stochastic, which further guarantees that $\lambda_1^l(t)$ is piecewise continuous and bounded, enabling the integrability of \eqref{eq_obsdeB}. Since it has been proven that $\tilde \G^l_{\mH}(t)$ always contains a directed spanning tree, Lemma~\ref{lm_obsever} can be used to ensure that the observer states of each normal agent will asymptotically converge to the same value and Zeno behavior will not occur, regardless of the behavior of Byzantine agents.

Moreover, in addition to being a sufficient condition, the (2$f$+1)-robustness of the graph $\G$ is also a necessary condition for this conclusion. If $\G$ is not (2$f$+1)-robust, it is possible to have $\C_{\mH} = \emptyset$, meaning that no convex hull made up of vertices only from normal agents can be found to enclose $\bar \W^l_i(t)$.

\begin{figure}[!t]
\centering
\includegraphics[width=0.33\linewidth]{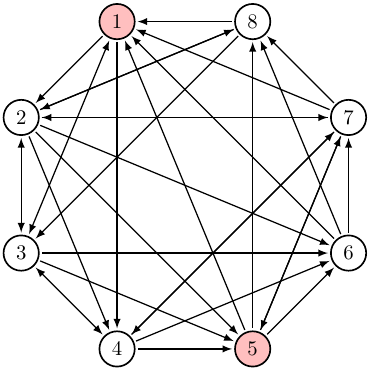}
\caption{\black{A 3-robust communication topology $\G$, where the red circles represent Byzantine agents and others are normal agents. The single-arrowed lines are utilized to depict directed communication links in $\G$, with the arrows pointing toward the agents receiving information.}}
\label{fig_topo}
\end{figure}

In conclusion, the observer state of each agent will asymptotically converge to the same value, and no Zeno behavior occurs. Next, the proof for achieving \eqref{eq_objcon} will be given. Choose the following Lyapunov function
\begin{equation}\label{eq_Ly}
V_{1,i}(t) = \frac{1}{2}s_{i}^T M_i s_{i} + \frac{1}{2F_i} \tilde \varphi_i^T \tilde \varphi_i + \frac{1}{2}e_{q\eta,i}^T e_{q\eta,i},
\end{equation}
where $\tilde \varphi_i = \hat \varphi_i - \varphi_i$ and $e_{q\eta,i}=q_i(t)-\eta_i(t)$. By using \eqref{eq_ELequation}, \eqref{eq_controller}, Property~\ref{p1} and $s_i = \dot q_i - v_i$, the derivative of \eqref{eq_Ly} can be found to be
\begin{align}\label{eq_dLyB1}
\dot V_{1,i}(t)
& { \;= } -k_i \|s_i\|^2 + (q_i-\eta_i)^T(s_i + v_i - \dot \eta_i) \notag\\
&\le -(k_i - \frac{1}{2c}) \|s_i\|^2 - (\mu_{2} - c)\|e_{q\eta,i}\|^2 \notag\\
&\quad + \frac{\alpha_6}{(t-t_0+\alpha_2)^{2 \alpha_3}},
\end{align}

\noindent
where $\alpha_6=\frac{\gamma_{S}^2\alpha_{v}^2}{2c}$ with $\alpha_{v}=d_\gamma \bar{\gamma}_J \bar{\gamma}_U \sqrt{n} + \mu_1 {\gamma}_{\tilde L} \alpha_1 \gamma_{-S}$. The second inequality in the above is obtained by using Young's inequality and the second equation of \eqref{eq_controller}. $c$ is a positive constant introduced by Young's inequality. From \eqref{eq_dLyB1}, it can be seen that there exists at least a constant $c$ such that $k_i >\frac{1}{2c}$ and $\mu_2 > c$ hold for any $k_i\mu_{2}>\frac{1}{2}$. In light of $\lim_{t \to \infty} \frac{\alpha_6}{(t-t_0+\alpha_2)^{2 \alpha_3}} = 0$, we have $\lim_{t \to \infty} \dot V_{1,i}(t) \le 0$ if $k_i\mu_{2}>\frac{1}{2}$, and $\dot V_{1,i}(t) = 0$ holds if and only if both $\|s_i\| = 0$ and $\|e_{q\eta,i}\| = 0$. Therefore, for any normal agent $i$, we have $\lim_{t \to \infty} \|q_{i}(t)-\eta_{i}(t)\|=0$ {and $\lim_{t \to \infty} \|\dot q_{i}(t)-v_{i}(t)\|=0$}. Since the observer state of each normal agent will asymptotically converge to the same value, i.e., $\lim_{t \to \infty} \|\eta_{j}(t)-\eta_{i}(t)\|=0$, we have $\lim_{t \to \infty} \| q_{j}(t)- q_{i}(t)\|=0$ for any pair of agents $i,j \in \mH$. Moreover, it follows from the definition of $v_i(t)$ that $\lim_{t \to \infty} v_i(t)=S\eta_{i}(t)$. Reviewing the definitions of $s_i(t)$ and $\dot \eta _i(t)$, it is straightforward to get that $\lim_{t \to \infty} \dot q_i(t)= v_{i}(t)$ and $\lim_{t \to \infty} \dot \eta_i(t)=S\eta_{i}(t)$, implying $\lim_{t \to \infty} \|\dot q_i(t)- \dot \eta_i(t)\|=0$. Combined with $\lim_{t \to \infty} \|\eta_{j}(t)-\eta_{i}(t)\|=0$, we can further obtain $\lim_{t \to \infty} \|\dot q_{j}(t)-\dot q_{i}(t)\|=0$.

Additionally, inspired by \eqref{eq_zeno1} and \eqref{eq_zeno2}, let us define $\bar{\alpha}_{z1}$ as a bounded constant greater than all possible values of $\alpha_z -\|S\| \alpha_1$. Then, by integrating \eqref{eq_dW2}, we obtain $\|\W_i(t)\|\le \|\W_i(t_0)\|+ \int \nolimits_{t_0}^{t} \frac{\bar{\alpha}_{z1}}{(\tau - t_0 +\alpha_2)^{\alpha_3}} d\tau$. It is evident that $\W_i(t)$ is bounded. Combined with its definition, this implies that $\|\eta_i(t)\|\le \gamma_S \|\W_i(t)\|$, i.e., $\eta_i(t)$ is also bounded. Similarly, considering that \eqref{eq_obs} can be rewritten as $\dot \eta_i(t) = S\eta_i(t) - \mu_1 e^{St}\sum \nolimits_{j\in \mH} \tilde a_{ij}(t) \circ (\hat {\W}_i(t) - \hat{\W}_j(t))$, we conclude that $\dot \eta_i(t)$ is also bounded. Furthermore, by integrating \eqref{eq_dLyB1}, we have that $s_i$ and $e_{q\eta,i}$ are bounded. Thus, $q_i(t)$ and $\dot q_i(t)$ are bounded.

So far, it can be concluded that the resilient consensus problem can be solved regardless of the behavior patterns of Byzantine agents.
\end{IEEEproof}

\section{Simulation Results}

In this section, we present some numerical examples to verify the effectiveness of the proposed algorithm. We consider a networked EL system comprising 8 two-link robotic arms, with each agent assumed to be modeled by \eqref{eq_ELequation}. Following \cite{models}, the system matrices of these agents are given as
\begin{align*}
M_i(q_i)&=\!\begin{bmatrix} \ell_{i1}+\ell_{i2}+2\ell_{i3}\cos q_{i2} & \ell_{i2}+\ell_{i3}\cos q_{i2} \\
\ell_{i2}+\ell_{i3}\cos q_{i2} & \ell_{i2} \end{bmatrix}\!\!, \\
C_i(q_i,\dot q_i)&=\!\begin{bmatrix} -\ell_{i3}\dot q_{i2}\sin q_{i2} & -\ell_{i3}(\dot q_{i1}+\dot q_{i2})\sin q_{i2} \\
\ell_{i3}\dot q_{i1}\sin q_{i2} & 0 \end{bmatrix}\!\!, \\
g_i(q_i)&=\!\begin{bmatrix} \ell_{i4}g\cos q_{i1}+\ell_{i5}g\cos (q_{i1}+q_{i2}) \\
\ell_{i5}g\cos (q_{i1}+q_{i2})\end{bmatrix}\!\!,
\end{align*}

\noindent
where $\tau_i=[\tau_{i}^1,\tau_{i}^2]^T$, $q_i = [q_{i}^1,q_{i}^2]^T$, and $\ell_i=[\ell_{i1},\ell_{i2},\ell_{i3},\ell_{i4},\ell_{i5}]^T$ are defined as the control input, the position state, and the system parameters of agent $i$, respectively. $g=9.8$ is the gravitational constant.

The system matrix of the observer of each agent is $S=[0,6;-1.5,0]^T$. The control gains of each agent are chosen as $\mu_1=5.9$, $\mu_2=2$, $k_i=80$, $\alpha_1=8$, $\alpha_2=3$, $\alpha_3=4$, and $F_i=0.6$. Let us define $\pi$ as pi, and the initial states and the physical parameters of agents are set as follows: $t_0=0$, $\ell_i = [0.64, 1.10, 0.08, 0.64, 0.32]^T$, $q_i(t_0) = [0.1\pi(i-1), -0.1\pi(i-11)]^T$, $\eta_1(t_0) = [-1.5, -0.5]^T$, $\eta_2(t_0) = [1, 0.5]^T$, $\eta_3(t_0) = [0, 0]^T$, $\eta_4(t_0) = [0.5, -2]^T$, $\eta_5(t_0) = [2,-1]^T$, $\eta_6(t_0) = [1.5, -0.5]^T$, $\eta_7(t_0) = [-1.5, -1]^T$, and $\eta_8(t_0) = [-2, -2]^T$. Besides, set $\dot q_i(t_0)=[0,0]^T,\;\forall i\in \V$.

The communication topology is shown in Fig.~\ref{fig_topo}, depicting a 3-robust digraph under 1-local Byzantine attacks, with agents 1 and 5 designated as Byzantine agents. The behavior patterns of Byzantine agents are given as follows.
Byzantine agent 5 simulates a fault agent with a fault on its control input. Considering that the derivative of the observer state of agent 5 is computed by \eqref{eq_obs} at time $t$ is $\dot \eta_5^c(t)$, but its practical available derivative of the observer state is designed as $\dot \eta_5(t)= \dot \eta_5^c(t)\circ[0.2, \sin(5t)+1]^T$. The state trajectory of Byzantine agent 5's observer is updated according to $\dot \eta_5(t)=[2\cos(4t), 0.4(\cos(8t)+\sin(6t))]^T$, regardless of the information it receives. Additionally, it will disconnect from its all out-neighbors when $t>t_0$. Byzantine agent 1 computes the derivative of the observer state by \eqref{eq_obs}, but it can arbitrarily design the information transmitted to its each out-neighbor. Let $\eta_{1}^{1\rightarrow j}(t)$ be the information transmitted from agent $1$ to its out-neighbor $j$, which is designed as
\begin{equation}\label{eq_attacks6}
\eta_{1}^{1\rightarrow j}(t)=\left\{\begin{array}{lll}
0.3j\eta_{1}(t), &\!j\in\{2\}, \\
\eta_{1}(t) + \alpha_{1j} + \beta_{1j}, &\!j\in\{3,4\},
\end{array} \right.
\end{equation}
in which, $\alpha_{1j}=[\alpha_{1j}^1,\alpha_{1j}^2]^T$ is the injected false data with $\alpha_{1j}^k=\min\{\sin(jt)|\eta_{j}^k|,|\eta_{1}^k(t_0)|\}$, $k\in\{1,2\}$, and $\beta_{1j}=[0.5\cos(t),0.05\cos(t)]^T$ is the measurement noise. And when $t>t_0$, it will disconnect from agent 4. Moreover, these Byzantine agents transmit their states to their out-neighbors at intervals of 0.001s. Obviously, such attacks violate the constraint that Byzantine agents cannot disconnect from their out-neighbors within continuous-time scenarios, rendering the theoretical analysis in \cite{vecbook,MSRE3,MSRE3n,MSRE4} ineffective.

%

\begin{figure*}[t!]
\centering
\subfigure[{\scriptsize Position state evolution}]{
    \includegraphics[width=0.225\textwidth]{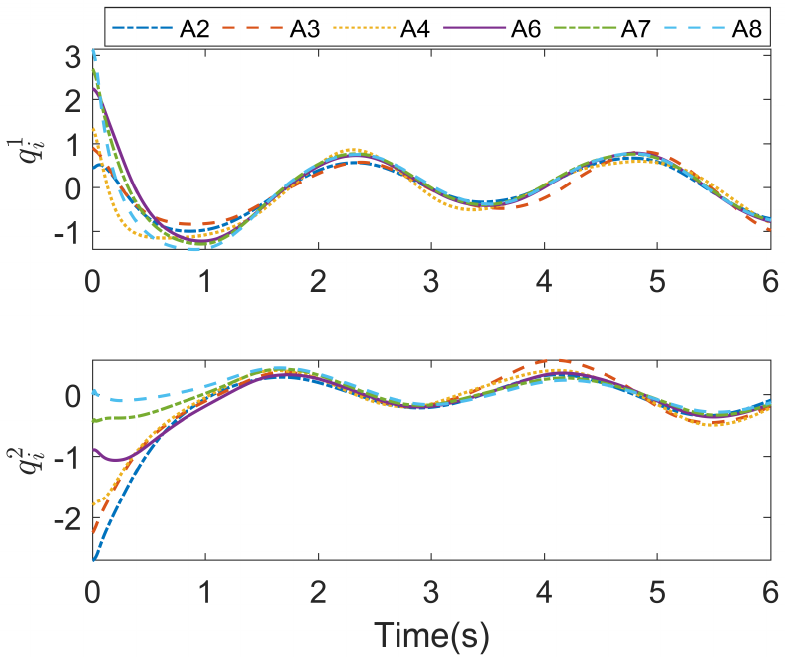}}\!\!~
\subfigure[{\scriptsize Velocity state evolution}]{
    \includegraphics[width=0.225\textwidth]{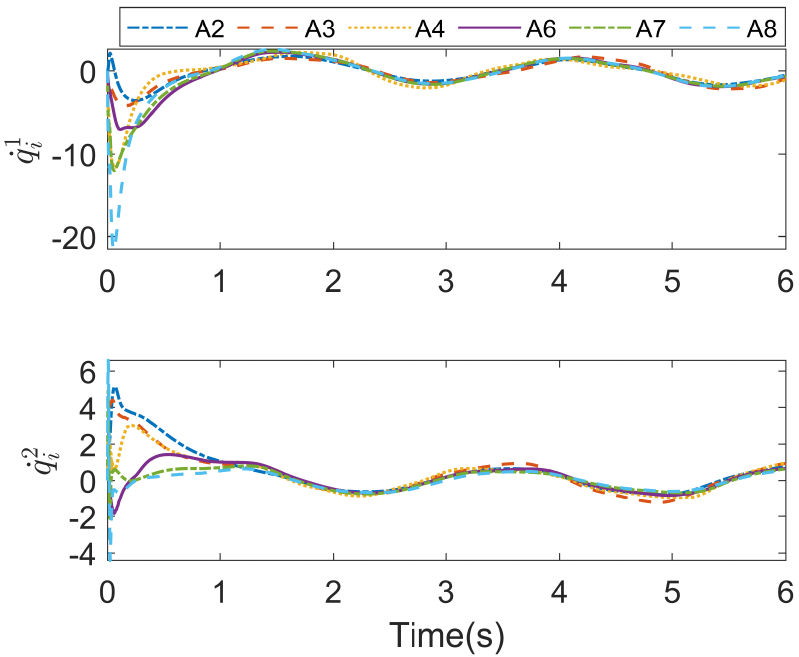}}\!\!~
\subfigure[{\scriptsize Auxiliary variable evolution}]{
    \includegraphics[width=0.225\textwidth]{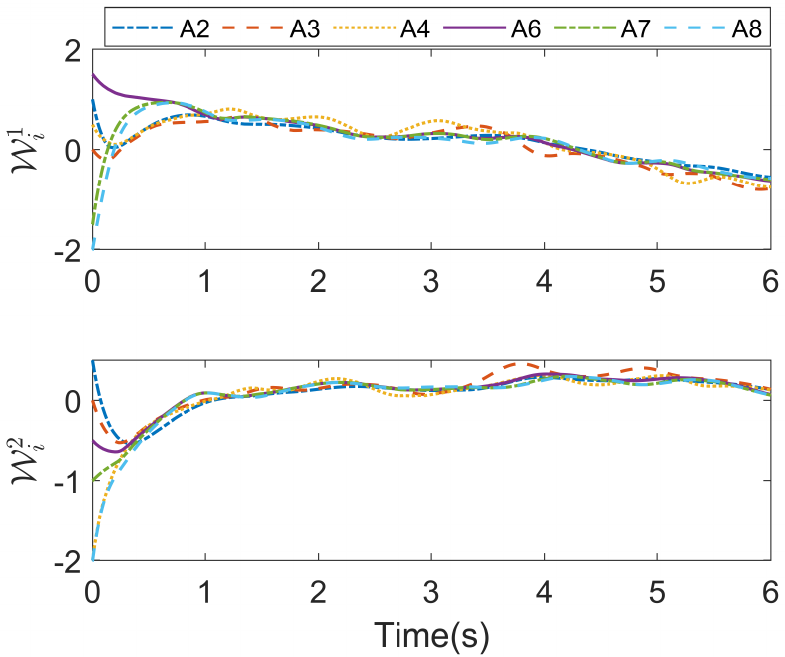}}\!\!~
\subfigure[{\scriptsize Triggering instants}]{
    \includegraphics[width=0.29\textwidth,height=23ex]{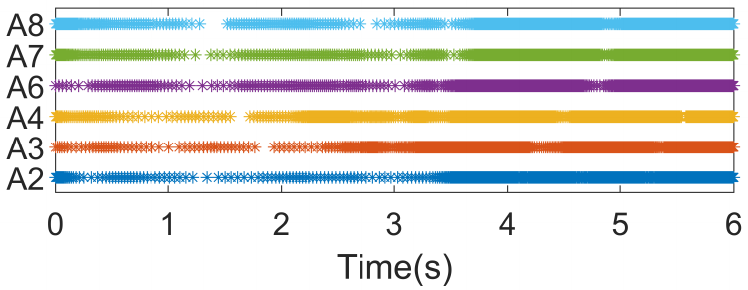}}
\caption{\black{\scriptsize Simulation results of each normal agent using Algorithm~\ref{alg-2} with $f=0$.}}
\label{qdWEno}
\end{figure*}

\begin{figure*}[t!]
\centering
\subfigure[{\scriptsize Position state evolution}]{
    \includegraphics[width=0.225\textwidth]{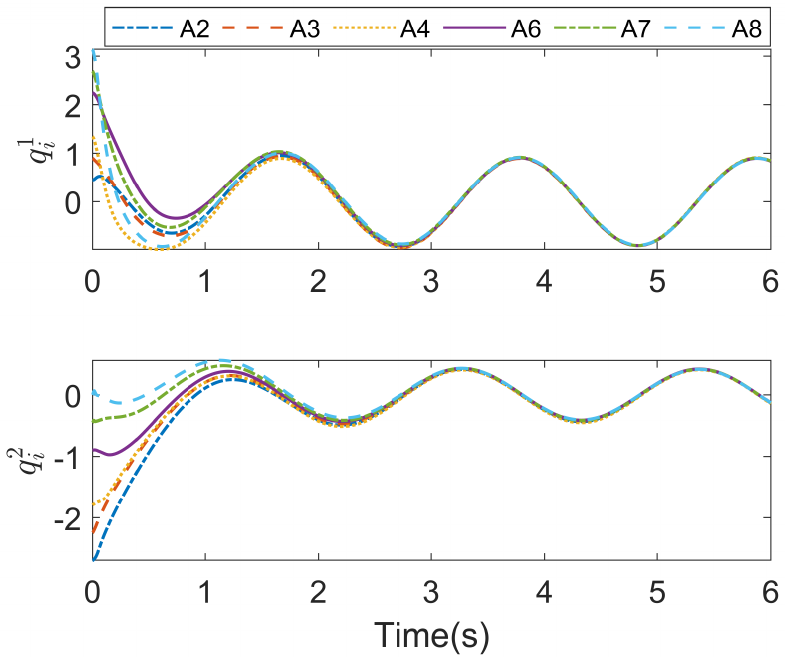}}\!\!~
\subfigure[{\scriptsize Velocity state evolution}]{
    \includegraphics[width=0.225\textwidth]{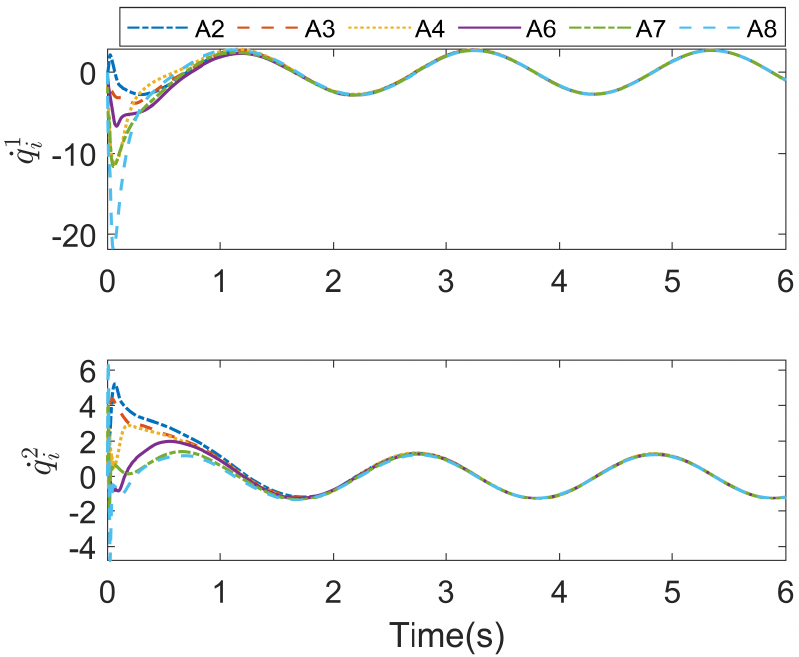}}\!\!~
\subfigure[{\scriptsize Auxiliary variable evolution}]{
    \includegraphics[width=0.225\textwidth]{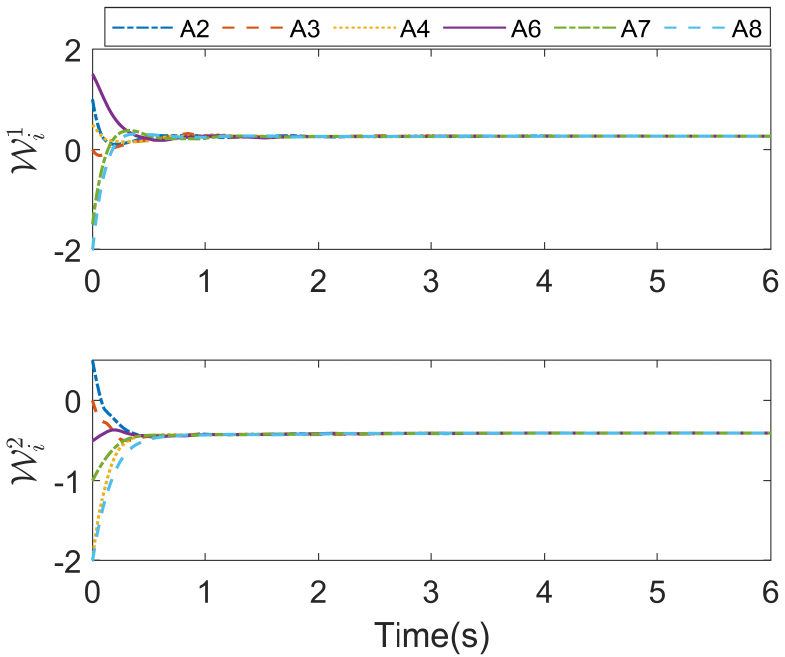}}\!\!~
\subfigure[{\scriptsize Triggering instants}]{
    \includegraphics[width=0.29\textwidth,height=23ex]{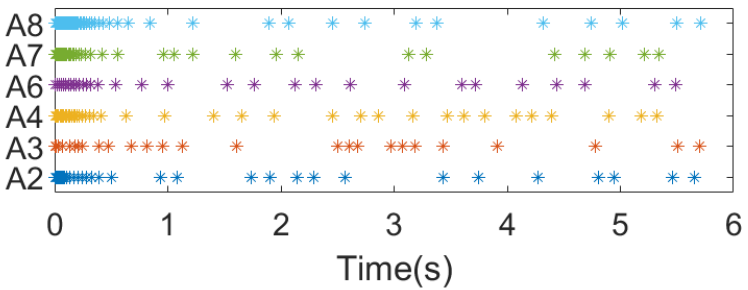}}
\caption{\black{\scriptsize Simulation results of each normal agent using Algorithm~\ref{alg-2} with the decision algorithm in \cite{vce4}.}}
\label{qdWEcmpr}
\end{figure*}

The simulation results are shown in {Figs.~\ref{qdWEno}--\black{\ref{qdWEcmpr}}} and Table~\ref{table1}, where the index $Ai$ is used to represent agent $i$ with $i \in \V$. Fig.~\ref{qdWEno}(a)\--(c) illustrate the scenarios that the AVBRD algorithm does not work, which is done by setting $f=0$ when applying Algorithm~\ref{alg-2}. From Fig.~\ref{qdWEno}(a)--(b), it can be observed that the designed Byzantine agents can easily disrupt the consensus of the networked EL system in that case. As a result, the communication of each agent is frequently triggered, as shown in Fig.~\ref{qdWEno}(d). Accordingly, as depicted in Fig.~\ref{qdWEno}(c), the designed auxiliary variables also fail to achieve consensus, meaning that the observer states do not converge either. These results demonstrate the effectiveness of the designed attacks.

%

\begin{figure*}[t!]
\centering
\subfigure[{\scriptsize Position state evolution}]{
    \includegraphics[width=0.225\textwidth]{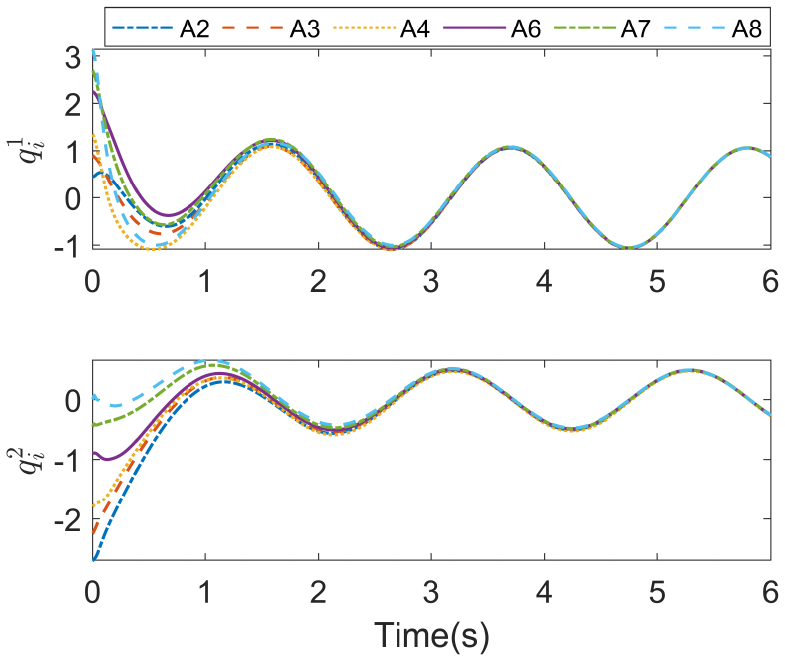}}\!\!~
\subfigure[{\scriptsize Velocity state evolution}]{
    \includegraphics[width=0.225\textwidth]{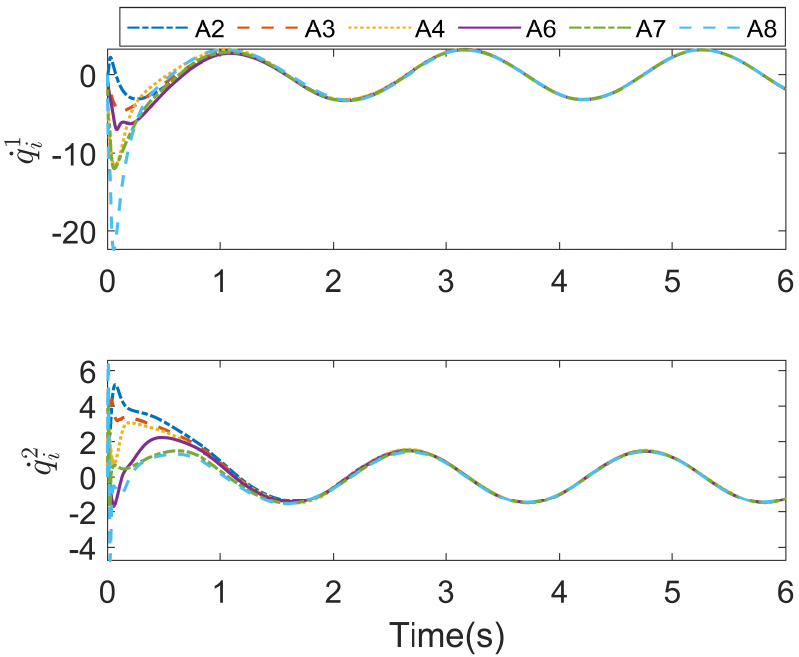}}\!\!~
\subfigure[{\scriptsize Auxiliary variable evolution}]{
    \includegraphics[width=0.225\textwidth]{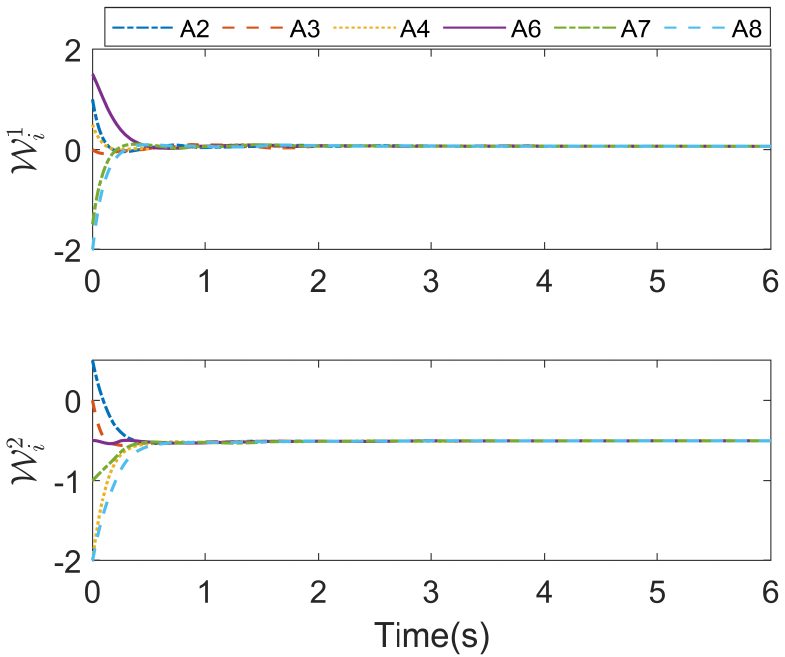}}\!\!~
\subfigure[{\scriptsize Triggering instants}]{
    \includegraphics[width=0.29\textwidth,height=23ex]{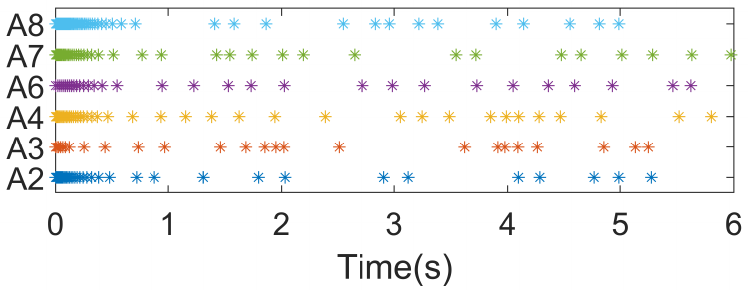}}
\caption{\black{\scriptsize Simulation results of each normal agent using Algorithm~\ref{alg-2} with $f=1$.}}
\label{qdWEmy}
\end{figure*}

%

Next, we set $f=1$ to verify the main results in Theorem~\ref{them1}, with the simulation results being shown in Fig.~\ref{qdWEmy}(a)--(c). From Fig.~\ref{qdWEmy}(c), it can be concluded that the proposed algorithm effectively ensures consensus among the designed auxiliary variables, as analyzed in Lemma~\ref{lm_obsever}. This means that the observer states of normal agents consequently achieve consensus, as validated in Table~\ref{table1}, where $e_{\eta_{ij},i}=\sum \nolimits_{j \in \N_i \bigcap \mH} (\eta_{i}-\eta_j)$. Moreover, as displayed in Fig.~\ref{qdWEmy}(a)--(b), the states of normal agents can achieve consensus and remain bounded. Additionally, from Fig.~\ref{qdWEmy}(d) and Table~\ref{table1}, it can be seen that the designed ET scheme can effectively reduce the communication resource consumption and ensure a minimum triggering interval. Altogether, these results confirm that the resilient consensus problem as described in Problem~\ref{df_obj1} is resolved.
\black{Furthermore, comparative simulations are conducted to verify the efficiency of Algorithm~\ref{alg-1} in reducing computational resources.
Specifically, in Step 2 of Algorithm~\ref{alg-2}, agent $i$ computes the centerpoint of $[ \hat {\W}_{j_1}(t), \hat {\W}_{j_2}(t), \cdots, \hat {\W}_{j_|\N_i|}(t)]$ as defined in \cite[Theorem~4.3]{vce4}, instead of executing Algorithm~\ref{alg-1}, to obtain $\bar{\W}_i(t)$. The results are exhibited in Fig.~\ref{qdWEcmpr}. Since $\bar{\W}_i(t)$ still lies within the convex hull of the corresponding states of agent $i$'s normal in-neighbors, resilient consensus can be achieved. However, executing the decision algorithm 26417 times takes a total of 3.749s, whereas Algorithm~\ref{alg-1} takes only 0.617 seconds for 24813 executions, showing that each execution of Algorithm~\ref{alg-1} uses less computational resources, as noted in Remark~\ref{rmk32}.
Moreover, even when Byzantine agents transmit information every 0.001s, the number of executions of Algorithm~\ref{alg-1} by any agent does not exceed 24813, which significantly reduces computational resources compared with executing it at every instant, thereby validating the statement in Remark~\ref{rmk31}.} So far, the effectiveness of the proposed resilient consensus algorithm is verified.

\begin{table}[!t]
\caption{Performance of the proposed event-triggered resilient consensus algorithm
\label{table1}}
\setlength\tabcolsep{4pt}
\renewcommand\arraystretch{1.2}
\centering
\begin{threeparttable}
\begin{tabular}{c|c|c|c|c|c|c}
   \hline
   {Agent index}           & A2      & A3      & A4      & A6      & A7      & A8     \\ \hline 
   {Triggering numbers}      & \black{36}      & \black{25}      & \black{40}      & \black{32}      & \black{41}      & \black{56}     \\ \hline
   {Settling time\tnote{1}}  & \black{1.051s}  & \black{1.961s}  & \black{1.428s}  & \black{1.693s}  & \black{0.572s}  & \black{0.632s} \\ \hline
   {\makecell{Minimum \\ triggering interval}}      & \black{0.007s}      & \black{0.013s}     & \black{0.009s}      & \black{0.017s}      & \black{0.006s}      & \black{0.005s}     \\ \hline
\end{tabular}
\begin{tablenotes}
\footnotesize
\item[1] Settling time: the time required for $\|e_{\eta_{ij},i}\|$ of agent $i$ to converge within 1\% of its initial value.  
\end{tablenotes}
\end{threeparttable}
\end{table}

\section{Conclusion}

This paper has first proposed a new resilient decision algorithm, named as the auxiliary-variable-based resilient decision (AVBRD) algorithm, to solve the resilient consensus problem of networked EL systems in the presence of Byzantine attacks.
The proposed algorithm can relax the constraints imposed on Byzantine agent behavior patterns in existing results within continuous-time scenarios by utilizing the minimal triggering interval and the open-loop estimation mechanism of the ET communication scheme.
The rigorous proofs and case studies have been made, verifying the effectiveness of the proposed algorithm. Future work will involve improving the proposed algorithm to reduce the requirement on graph robustness.

\section*{Appendix}

\subsection{Proof of Lemma~\ref{lm_neighbours}}

\begin{IEEEproof}[\hspace{-1.0em}Proof]
Recalling Definition~\ref{dfrobust}, if $S_1=\{i\}$ and $S_2= \V \backslash \{i\}$, then the $r$-reachable set must be $S_1$, as there is only one agent outside of $S_2$. Therefore, according to Definition~\ref{dfreachable}, we have $|\N_i\backslash \{i\}| \ge r$, which further implies $|\N_i| \ge r$. This completes the proof.
\end{IEEEproof}

\subsection{Proof of Lemma~\ref{lm_obsever2lm}}

\begin{IEEEproof}[\hspace{-1.0em}Proof]
First, let us discuss the case where $f \neq 0$. In this scenario, it follows from Lemma~3.3 that $|\N_i| \ge 2f+1$ and $|\N_i \bigcap \F| \le f$.
Therefore, by constructing convex hulls using the auxiliary variable values from each subset of $|\N_i| - f$ agents in $\N_i$, and taking the intersection of all possible convex hulls, one can obtain a convex hull $\C^l_{\mH,i}$ contained within the convex hull formed by the auxiliary variable values of all normal agents on dimension $l$.
Considering that for a set of scalars, the convex hull is either a point or a line segment defined by its maximum and minimum values, one has $\C^l_{\mH,i}= \C^l_{m,i} \bigcap \C^l_{M,i}$, where $\C^l_{m,i}$ and $\C^l_{M,i}$ are the convex hulls formed by the smallest and largest $|\N_i| - f$ auxiliary variables, respectively, sharing $|\N_i|-2f$ vertices. Thus, $\C^l_{\mH,i} \neq \emptyset$ and $\C^l_{\mH,i}$ contains these shared vertices when $|\N_i|\ge 2f+1$. Further, if $\C^l_{\mH,i}$ has $k$ vertices from Byzantine agents, then both $\C^l_{m,i}\backslash\C^l_{\mH,i}$ and $ \C^l_{M,i}\backslash\C^l_{\mH,i}$ have at least $k$ vertices from normal agents. These vertices can form $k$ convex hulls, such that the vertices from Byzantine agents in $\C^l_{\mH,i}$ can uniquely lie inside or on the vertices of one of these convex hulls. Thus, any point in $\C^l_{\mH,i}$ lies in a convex hull formed by at least $|\N_i|-2f$ vertices from normal agents.
Moreover, since the vertices shared by $\C^l_{m,i}$ and $ \C^l_{M,i}$ have the most central state values, $\bar{\W}^l_i(t)$ is actually the mean of the maximum and minimum vertices of $\C^l_{\mH,i}$, implying $\bar{\W}^l_i(t) \in \C^l_{\mH,i}$. The aforementioned conclusion thus holds for $\bar{\W}^l_i(t)$. Based on a similar process, this conclusion is obvious when $f=0$. So far, the proof is complete.
\end{IEEEproof}

\subsection{Proof of Lemma~\ref{lm_obsever}}

\begin{IEEEproof}[\hspace{-1.0em}Proof]
It is evident that the dynamics described in \eqref{eq_lmobseq} are decoupled across dimensions. Utilizing this, we can derive $\dot \W_i^l(t) = - \mu_1 \sum \nolimits_{j=1}^{N} \tilde a_{ij}^l(t)(\hat {\W}^l_i(t) - \hat{\W}^l_j(t))$, where $\W^l_i(t)$ and $\hat {\W}^l_i(t)$ are the $l$-th dimension of $\W_i(t)$ and $\hat {\W}_i(t)$, respectively. Thus, \eqref{eq_lmobseq} can be written as the following compact form
\begin{equation}\label{eq_obseqcom}
\dot \W^l(t)= -\mu_1 \tilde L^l \W^l(t) -\mu_1 \tilde L^l(t) e_{\W}^l(t),
\end{equation}
where $\W^l(t)=[ {\W}_{1}^l(t), {\W}_{2}^l(t), \cdots, {\W}_{N}^l(t)]^T$, $e_{\W}^l(t)=[ e_{\W 1}^l(t), e_{\W 2}^l(t), \cdots, e_{\W N}^l(t)]^T$ with $e_{\W i}^l(t)$ is $l$-th dimension of $e_{\W i}(t)$ and $e_{\W i}(t)=e^{-St} e_{\eta i}(t)$, and $\tilde L^l(t)$ is the Laplacian matrix of $\tilde \G^l(t)$.
Let ${J_{\tilde L^l}(t)}\in \mathbb{R}^{N\times N}$ be the Jordan canonical form of ${\tilde L^l(t)}$, and $U^l(t) \in \mathbb{R}^{N\times N}$ be a nonsingular matrix such that $(U^l(t))^{-1}{{\tilde L^l}(t)}U^l(t)={J_{\tilde L^l}(t)}$. Let {$\epsilon^l(t)=(U^l)^{-1}\W^l(t)$}. Then \eqref{eq_obseqcom} can be rewritten as
\begin{equation}\label{eq_obsre}
{\dot \epsilon^l(t)= - \mu_1 J_{\tilde L^l}(t) \epsilon^l(t) - \mu_1 J_{\tilde L^l}(t) e_{\epsilon}^l(t),}
\end{equation}
where $e_{\epsilon}^l(t)=(U^l(t))^{-1}e_{\W}^l(t)$. Lemma~\ref{lm_seclm1} indicates that ${\tilde L^l(t)}$ has only one zero eigenvalue. Therefore, ${J_{\tilde L^l}(t)}$ can be rearranged as $J_{\tilde L^l}(t)=\block \, \diag(0,J^l_{2}(t))$, where $J^l_2(t) \in \mathbb{R}^{(N-1) \times (N-1)}$. Let $\epsilon^l_1(t)$ denote the first $1$ rows of $\epsilon^l(t)$, and $\epsilon^l_2(t)$ is the last $N-1$ rows of $\epsilon^l(t)$. Similarly, define $\dot \epsilon^l_1(t)$ and $\dot \epsilon^l_2(t)$ as the first $1$ rows and the last $N-1$ rows of $\dot \epsilon^l(t)$, respectively. $e^l_{\epsilon 2}(t)$ is given as the last $N-1$ rows of $e^l_{\epsilon}(t)$. Then \eqref{eq_obsre} is re-expressible as
\begin{subequations}\label{eq_obsde}
\begin{align}
\dot \epsilon^l_1(t) &=0, \label{eq_obsdeA}\\
\dot \epsilon^l_2(t) &=- \mu_1 J^l_2(t) \epsilon^l_2(t) - \mu_1 J^l_2(t) e^l_{\epsilon 2}(t). \label{eq_obsdeB}
\end{align}
\end{subequations}

Obviously, system \eqref{eq_obsdeA} is stable. Next, the analysis on the convergence of system \eqref{eq_obsdeB} will be given. From \eqref{eq_obsdeB}, by solving the differential equation, one can obtain
\begin{align}\label{eq_diff2}
\| \epsilon^l_2(t)\| &\le e^{-\lambda^l_1(t)(t-t_0)} \|\epsilon^l_2(t_0)\| \notag \\
&\quad + \gamma^l_J(t) \!\!\int \nolimits_{t_0}^{t} e^{-\lambda^l_1(t)(t-{\tau})} \|e^l_{\epsilon 2} (\tau)\| d\tau,
\end{align}

\noindent
where $\gamma^l_J(t){=}\mu_1 \|J^l_2(t)\|$, and $\lambda_1^l(t){=}\mu_1 \lambda_{J^l_2}(t)$ with $\lambda_{J^l_2}(t)$ being the smallest eigenvalue of ${J^l_2}(t)$. Since $\|e^l_{\epsilon 2}(\tau)\|{\le}\|e_{\epsilon 2}(\tau)\|$ and the eigenvalues of ${-}S$ are all considered with zero real parts, a positive constant $\gamma_{{-}S}$ can always be found such that $\| e^{-S(t-t_0)} \| \le \gamma_{{-}S} $ for any $t >t_0$.
Recalling the ET scheme \eqref{eq_event}, one has
\begin{align}\label{eq_diff4}
\| \epsilon^l_2(t)\| &\le e^{-\lambda^l_1(t)(t-{t_0})} \|\epsilon^l_2(t_0)\| \notag\\
&\quad + \int \nolimits_{t_0}^{t} e^{-\lambda^l_1(t) ({t-\tau})} \frac{\gamma^l_2(t)}{(\tau - t_0 +\alpha_2)^{\alpha_3}} d\tau,
\end{align}

\noindent
where $\gamma^l_2(t) = \gamma^l_J(t) \gamma_{{-}S} \alpha_1 \|(U^l(t))^{{-}1}\|$. Let $\rho_2^l(t) = (t - t_0 +\alpha_2)^{\alpha_3} \rho^l_1(t)$ with $\rho^l_1(t_0){=} \|\epsilon^l_2(t_0)\|$. The right-hand side of \eqref{eq_diff4} can be obtained by solving the following differential equation
\begin{equation}\label{eq_diff6}
\dot \rho^l_2(t) = \left(\frac{\alpha_3}{t - t_0 +\alpha_2} - \lambda^l_1(t)\right) \rho^l_2(t) + \gamma^l_2(t).
\end{equation}

According to Lemma~\ref{lm_seclm2}, it can be concluded that there must exist $\gamma_\omega>0$ for any $\gamma^l_1(t) = \frac{\int \nolimits_{t_0}^{t}\lambda^l_1(\tau)d\tau}{2(t-t_0)}$ such that
\begin{equation}\label{eq_diff7}
\int \nolimits_{\tau}^{t} \frac{\alpha_3}{s - t_0 +\alpha_2} ds\le \frac{\gamma^l_1(t)}{2}(t-\tau) + \gamma_\omega.
\end{equation}

By using \eqref{eq_diff6} and \eqref{eq_diff7}, one obtains
\begin{align}\label{eq_diff8}
\rho_2^l(t) &{\le} e^{\gamma_\omega} \rho^l_2(t_0) e^{-\frac{\gamma^l_1(t)}{2}(t - t_0)} + \bar \gamma_2 e^{\gamma_\omega} \int \nolimits_{t_0}^{t} e^{-\frac{\gamma^l_1(t)}{2}(t- \tau)} d\tau \notag\\
&{\le} e^{\gamma_\omega} \rho^l_2(t_0) + \bar \gamma_2 e^{\gamma_\omega} \frac{2}{\underline{\gamma}_1},
\end{align}

\noindent
where $\bar \gamma_2$ is the upper bound of $\gamma^l_2(t)$ and $\underline{\gamma}_1$ is the positive lower bound of $\gamma^l_1(t)$. The existence of $\bar \gamma_2$ and $\underline{\gamma}_1$ is ensured by the invertibility of $U^l(t)$ and the directed spanning tree in $\tilde \G^l(t)$.
Let $d^l_\gamma =e^{\gamma_\omega} \rho^l_2(t_0) + \gamma_2 e^{\gamma_\omega} \frac{2}{\underline{\gamma}_1}$. From \eqref{eq_diff4} and \eqref{eq_diff6}, it can be seen that
\begin{equation}\label{eq_diff9}
\| \epsilon_2^l(t) \| \le \frac{d^l_\gamma}{(t - t_0 +\alpha_2)^{\alpha_3}},
\end{equation}
which leads to $\lim \nolimits_{t \to \infty} \| \epsilon_2(t) \|\,{=}\,0$ for any $\alpha_2\,{>}\,1$ and $\alpha_3\,{>}\,1$. From the above discussion, it can be obtained that, $\lim \nolimits_{t \to \infty} \eta^l(t)\,{=}\,U^l\epsilon_3^l(t)$, where $\epsilon^l_3(t) \,{=}\, [(\epsilon^l_1)^T(t), \mathbf{0}_{N-1}^T]^T$ and $\eta^l(t)\,{=}\,[ {\eta}_{1}^l(t), {\eta}_{2}^l(t), {\cdots}, {\eta}_{N}^l(t)]^T$. Let $\mathbf{1}_{N}$ be the eigenvector associated to the zero eigenvalue. Therefore, we have $\lim \nolimits_{t \to \infty} \W^l_i(t)= \mathbf{1}_{N}^T\epsilon^l_1(t)$, for any agent $i$. Thus, the auxiliary variables of each agent can asymptotically converge to the same value, and there holds
\begin{equation}\label{eq_diff10}
\sum \limits_{l=1}^n \|\tilde L^l(t)\W^l(t)\|^2 \le \bigg(\frac{d_\gamma \bar \gamma_U \bar \gamma_J \sqrt{n}}{\mu_1 (t {-} t_0 {+}\alpha_2)^{\alpha_3}}\bigg)^2 ,
\end{equation}
where $d_\gamma=\max_l\{d_\gamma^l\}$, and $\bar \gamma_J$ and $\bar \gamma_U$ are the upper bound of all possible $\gamma_J^l(t)$ and $\|\gamma_U^l(t)\|$, respectively. Furthermore, from the definition of $\W_i$, one can get
$\|\eta_i(t)-\eta_j(t)\|\le \gamma_S\|\W_i(t)-\W_j(t)\|$ with $ \|e^{S(t-t_0)}\| \le \gamma_S$, where the existence of $\gamma_S$ is guaranteed by the fact that all the eigenvalues of $S$ have zero real parts. So far, the analysis of the asymptotic convergence of $\eta_i(t)$ is completed.

Next, the proof about the exclusion of the Zeno behavior under the ET scheme \eqref{eq_event} will be given. From \eqref{eq_obs}, \eqref{eq_hateta}, and \eqref{eq_lmobseq}, one has
\begin{equation}\label{eq_zeno1}
\dot {e}_{\eta i}(t) = S{e}_{\eta i}(t) + \mu_1 \W_{\tilde L^l i}(t) + \mu_1 e_{\tilde L^l i}(t),
\end{equation}
where $\W_{\tilde L^l i}(t){=}[\tilde L^1_i(t) \W^1(t), \tilde L^2_i(t) \W^2(t), \cdots\!, \tilde L^n_i(t) \W^n(t) ]^T$ and $e_{\tilde L^l i}(t)=[\tilde L^1_i(t) {e}_{\W }^1(t), \tilde L^2_i(t) {e}_{\W }^2(t) , \cdots,\tilde L^n_i(t) {e}_{\W }^n(t) ]^T$ with $\tilde L^l_i(t)$ being the $i$th row of $\tilde L^l(t)$. From \eqref{eq_event}, \eqref{eq_diff9}, \eqref{eq_diff10} and \eqref{eq_zeno1}, the upper bound of $D^+ \|{e}_{\eta i}(t)\|$ can be obtained as
\begin{equation}\label{eq_zeno2}
D^+ \|{e}_{\eta i}(t)\| {=} \frac{{e}_{\eta i}^{T}(t){\dot e}_{\eta i}(t)}{\|{e}_{\eta i}(t)\|} {\le} \|{\dot e}_{\eta i}(t)\| {\le} \frac{\alpha_z}{(t{-}t_0+\alpha_2)^{\alpha_3}},
\end{equation}
in which, $\alpha_z=\|S\|\alpha_1+ d_\gamma \bar \gamma_J \bar \gamma_U \sqrt{n} + \mu_1 \gamma_{\tilde L}\alpha_1 \gamma_{-S}$, $\gamma_{\tilde L}$ is the upper bound of all possible $\|\tilde L^l(t)\|$, and the fact that ${\dot e}_{\eta i}(t) = S { {\hat \eta}}_ i(t) -{\dot \eta}_ i(t)$ is used. Then one has
\begin{align}\label{eq_zeno3}
\|{e}_{\eta i}(t)\| &\le \|{e}_{\eta i}(t_k^i)\| +\int \nolimits_{t_k^i}^{t} \frac{\alpha_z}{(s-t_0+\alpha_2)^{\alpha_3}} ds \notag\\
&= \alpha_4 \!\bigg[\frac{1}{(t_k^i{-}t_0{+}\alpha_2)^{\alpha_3-1}} {-}\frac{1}{(t{-}t_0{+}\alpha_2)^{\alpha_3-1}}\bigg]\!,
\end{align}

\noindent
where $\alpha_4=\frac{\alpha_z}{\alpha_3-1}$. According to \eqref{eq_event}, the lower bound of $t_{k+1}^i$ can be obtained by solving the following inequality
\begin{align}\label{eq_zeno4}
\|{e}_{\eta i}(t_{k+1}^i)\| &= \frac{\alpha_1}{(t_{k+1}^i-t_0+\alpha_2)^{\alpha_3}} \notag\\
\le \alpha_4 & \bigg[\frac{1}{(t_k^i{-}t_0{+}\alpha_2)^{\alpha_3-1}} {-} \frac{1}{(t_{k+1}^i{-}t_0{+}\alpha_2)^{\alpha_3-1}}\bigg]\!.
\end{align}

By using the equivalent transformation and the fact that $t_{k+1}^i > t_{k}^i \ge t_{0}$, we can derive from \eqref{eq_zeno4} that
\begin{equation}\label{eq_zeno5}
\frac{\alpha_1}{\alpha_4(t_{k+1}^i-t_0+\alpha_2) } +1 \le \bigg(1 + \frac{t_{k+1}^i-t_k^i}{t_{k}^i-t_0+\alpha_2}\bigg)^{\alpha_3}.
\end{equation}

Define $\tilde{t}_{k}^i=t_{k}^i-t_0+\alpha_2$. According to \eqref{eq_zeno5} and $\alpha_3>1$, the lower bound of the triggering interval can be described as
\begingroup
\allowdisplaybreaks
\begin{align}\label{eq_intervals}
t_{k+1}^i - t_{k}^i &\ge \bigg(\bigg(\frac{\alpha_1}{\alpha_4 \tilde{t}_{k+1}^i} +1\bigg)^{\frac{1}{\lceil \alpha_3 \rceil}} -1 \bigg)\tilde{t}_{k}^i\notag\\
&= \frac{\alpha_1\tilde{t}_{k}^i }{\alpha_4 \tilde{t}_{k+1}^i \sum\limits_{\iota=0}^{{\lceil \alpha_3 \rceil}-1}{\Big(\frac{\alpha_1}{\alpha_4 \tilde{t}_{k+1}^i}+1\Big)^{\frac{\iota}{\lceil \alpha_3 \rceil}} }}\notag\\
&\ge \frac{\alpha_1\tilde{t}_{k}^i}{\alpha_4 \tilde{t}_{k+1}^i \sum\limits_{\iota=0}^{{\lceil \alpha_3 \rceil}-1}{\Big(\frac{\alpha_1}{\alpha_4\alpha_2}+1\Big)^{\frac{\iota}{\lceil \alpha_3 \rceil}} }},
\end{align}
\endgroup

\noindent
where $\tilde{t}_{k+1}^i=t_{k+1}^i-t_0+\alpha_2$, the equality is obtained by applying the formula $\varpi_1^\kappa - \varpi_2^\kappa = (\varpi_1-\varpi_2) \sum\nolimits_{\iota=0}^{\kappa-1}{\varpi_1^{\kappa-1-\iota}\varpi_2^\iota}$ with $\varpi_1=1$, $\varpi_2=\left(\frac{\alpha_1}{\alpha_4 \tilde{t}_{k+1}^i} +1\right)^{\frac{1}{\lceil \alpha_3 \rceil}}$ and $\kappa=\lceil \alpha_3 \rceil$.

Let $\Delta t = t_{k+1}^i - t_{k}^i$ and $\alpha_5 = \sum\nolimits_{\iota=0}^{{\lceil \alpha_3 \rceil}-1}{\left(\frac{\alpha_1}{\alpha_4\alpha_2}+1\right)^{\frac{\iota}{\lceil \alpha_3 \rceil}}}$. Then one has
\begin{equation}\label{eq_intervals2}
\Delta t \ge \frac{\alpha_1}{\alpha_4 \alpha_5} \cdot \frac{1}{ \frac{\Delta t}{\tilde{t}_{k}^i} + 1} \ge \frac{\alpha_1}{\alpha_4 \alpha_5} \cdot \frac{\alpha_2}{ {\Delta t} + \alpha_2}.
\end{equation}
From \eqref{eq_intervals2}, one can obtain
\begin{equation}\label{eq_intervals3}
\Delta t (\Delta t + \alpha_2)\ge \frac{\alpha_1 \alpha_2}{\alpha_4 \alpha_5}>0,
\end{equation}
from which it can be concluded that there is a positive constant $\varrho$ such that $t_{k+1}^i-t_k^i \ge \varrho$. Thereby, as long as $\W_i(t_0)$ is kept bounded for any agent $i$, $\tilde \G^l(t)$ satisfies the conditions described in Lemma~\ref{lm_obsever}, and no Zeno behavior occurs at the initial instant, one has that the Zeno behavior can be eradicated and the convergence of $\eta_i(t)$ can be ensured, regardless of the specific switching instants and frequencies of $\tilde \G^l(t)$. So far, the proof about Lemma~\ref{lm_obsever} is complete.
\end{IEEEproof}



\end{document}